\documentclass[12pt,review]{elsarticle}
\usepackage{ctex}
\usepackage{epstopdf}
\usepackage[caption=false]{subfig}
\usepackage{marvosym}
\usepackage{amsthm}
\usepackage{amssymb}
\usepackage{stackrel}
\usepackage{xcolor}
\usepackage{amsmath,amscd}

\newtheorem{dl}{Theorem}[section]
\newtheorem{tl}[dl]{Corollary}
\newtheorem{yl}[dl]{Lemma}
\newtheorem{dy}[dl]{Definition}
\newtheorem{lz}[dl]{Example}

\newtheorem{xz}[dl]{Proposition}

\newtheorem{remark}[dl]{Remark}

\usepackage[colorlinks=true,linkcolor=green,filecolor=brown,citecolor=green]{hyperref}
\usepackage[active]{srcltx}

\newproof{pot333}{Proof of Theorem \ref{firsteq}}
\newproof{pot444}{Proof of Theorem \ref{secondeq-secondnew}}
\newproof{pot777}{Proof of Theorem \ref{secondeq-third-new}}

\newproof{pot555}{Applications of Theorem \ref{firsteq}}
\newproof{pot666}{Applications of Theorem \ref{secondeq-secondnew}}
\newproof{pot888}{Applications of Theorem \ref{secondeq-third-new}}
\newcommand{\poq}[2]{(#1;q)_{#2}}

\def\qed{\hfill \rule{4pt}{7pt}}
\def\pf{\noindent {\it Proof.} }

\title{Bilinear generating functions of the multivariable Al-Salam-Carlitz polynomials  and applications}
\author{Qi Chen $^{a,}$\fnref{fn1}}
\fntext[fn1]{ E-mail address: qchen2025@stu.suda.edu.cn}
\address[P.R.China]{Department of Mathematics, Soochow University, Suzhou 215006, P. R. China}
\author{Xinrong Ma $^{a,}$\fnref{fn2}}
\fntext[fn2]{ E-mail address: xrma@suda.edu.cn}
\author{Jin Wang $^{b,}$\fnref{fn3}}
\fntext[fn3]{Corresponding author. E-mail address: jinwang@zjnu.edu.cn}
\address[P.R.China]{School of Mathematical  Sciences, Zhejiang Normal University, Jinhua 321004, P. R. China}
\begin{document}
\begin{abstract}In this paper, by the method of  comparing coefficients, we  establish a new generating function of multivariable Al-Salam-Carlitz polynomials which contains both Rogers' and Bowman's symmetirc expansion formulas and the classical $q$-Mehler formula  as special cases.  Some new $q$-series identities related to the multivariable Al-Salam-Carlitz polynomials are also presented.
\end{abstract}
\maketitle
\section{Introduction}
Throughout this paper, we will follow the notation and
terminology in \cite{10}.  The $q$-shifted factorials with $|q|<1$ and the variable $a$ are defined by
\begin{align*}
(a;q)_\infty
:=\prod_{n=0}^{\infty}(1-aq^n),\
(a;q)_n:=\frac{(a;q)_\infty}{(aq^n;q)_\infty}
\end{align*}
for all $n\in\mathbb{Z}, \mathbb{Z}$ is the set of integers.  We also adopt the notation
\[(a_1,a_2,\cdots,a_m;q)_n :=(a_1;q)_n (a_2;q)_n\cdots (a_m;q)_n.\]
Meanwhile, the $q$-binomial coefficient is defined by
\begin{align*}
\genfrac{[}{]}{0pt}{}{n}{k}_{q}:=\left\{
  \begin{array}{cc}
 \displaystyle  \frac{\poq{q}{n}}{\poq{q}{k}\poq{q}{n-k}},&\mbox{if}~~0\leq k\leq n;\\
0,&\mbox{otherwise}.
   \end{array}
\right.
\end{align*}
The  basic hypergeometric series with the base $q$  and the argument $z$ is
defined by
\begin{align*}{}_{r}\phi_s\left[\begin{matrix}a_1,&a_2,&\ldots,&a_{r}  \\  &b_1,&\ldots,&b_s \end{matrix};q,z \right]
	:=\sum_{n=0}^\infty
	\frac{(a_1,a_2,\ldots,a_{r};q)_n}{(q,b_1,\ldots,b_s;q)_n}(\tau(n))^{s+1-r} z^n,\end{align*}
 where $\tau(n):=(-1)^nq^{n(n-1)/2}, s+1\geq r$.

As a kind of  important orthogonal polynomials (cf. \cite{[2], [4], [3]}), the Al-Salam-Carlitz polynomials $\Phi_n^{(\alpha)}(x | q)$ in \cite[(1.11)]{111} are  defined as
\begin{align*}
\Phi_n^{(\alpha)}(x | q)&:=\sum_{k=0}^n\genfrac{[}{]}{0pt}{}{n}{k}_{q}(\alpha;q)_k x^k.
\end{align*}
It is to be noted that the authors of the paper \cite{141414}  introduced the homogeneous  form of the Al-Salam-Carlitz polynomials
\begin{align}
\Phi_n^{(\alpha)}(x, y | q):=\sum_{k=0}^n\genfrac{[}{]}{0pt}{}{n}{k}_{q}(\alpha ; q)_k x^k y^{n-k}.\label{hanpoly}
\end{align}
It is clear that the following basic relations hold
\begin{align*}
\Phi_n^{(\alpha)}(x, y | q)&=y^n \Phi_n^{(\alpha)}(x / y | q);\qquad
\Phi_n^{(\alpha)}(x, 1 | q)=\Phi_n^{(\alpha)}(x | q).
\end{align*}
Observe that when $\alpha=0$ in \eqref{hanpoly}, these polynomials coincide with the classical  homogeneous Rogers-Szeg\"{o} polynomials
\begin{align}
h_n(x, y | q)=\sum_{k=0}^n\genfrac{[}{]}{0pt}{}{n}{k}_{q} x^k y^{n-k}.\label{h-melher}
\end{align}
To our knowledge, one of the important results about the homogeneous Al-Salam-Carlitz polynomials is the following bilinear generating function (in short, GF) due to Al-Salam and Carlitz (cf. \cite{111, 2015}), which is also called the $q$-Mehler formula for the homogeneous Al-Salam-Carlitz  polynomials.
\begin{dl}[Cf. \mbox{\cite[(1.9)]{111} and \cite[Thm. 3.2]{2015}}]\label{mainthm-thm} For $\max\{|x t u|,|x t v|,|y t u|,$ $|y t v|\}<1$, we have
\begin{align}
\sum_{n=0}^{\infty} \Phi_n^{(\alpha)}(x, y | q) \Phi_n^{(\beta)}(u, v | q) \frac{t^n}{(q ; q)_n}=\frac{(\alpha x t v, \beta y t u ; q)_{\infty}}{(x t v, y t v, y t u ; q)_{\infty}}{}_3 \phi_2\left[\begin{matrix}
    \alpha, \beta, y t v  \\
\alpha x t v, \beta y t u
\end{matrix}; q, x t u\right].\label{mainthm}
\end{align}
\end{dl}
So far, there have been a few different approaches to this formula \cite{111,555,  [3],101010}. For instance, Ismail and Stanton   \cite[p. 1042, (5.17)]{[4]} found the special result when $\alpha=\beta=0$. Liu \cite[Thm. 3.2]{2015} showed that it can be recognized as a solution of the certain $q$-partial differential equation.  Srivastava, Cao, and Jain \cite{363636, 141414} investigated   multi-linear generating functions of  the $q$-Hahn and $q$-Hermite polynomials.

In our previous paper \cite{chenwang}, we established the following generating function of the Al-Salam-Carlitz polynomials.
\begin{dl}[Cf.
{\rm \cite[Thm. 2.1]{chenwang}}]\label{temper} For any integer $m\geq 0$, it holds
\begin{align}
\sum_{n=0}^{\infty} \Phi_{n+m}^{(\alpha)}(x, y | q) \frac{t^{n}}{(q ; q)_n}=\frac{(\alpha x t; q)_{\infty}}{(x t, y t; q)_{\infty}}\sum_{n=0}^m \genfrac{[}{]}{0pt}{}{m}{n}_{q}
\frac{(\alpha,  y t;q)_n}{(\alpha x t;q)_n}x^ny^{m-n}.\label{mainthm-I}
\end{align}
\end{dl}

 As one of the main results in \cite{chenwang}, we showed that the GF  \eqref{mainthm-I} also implies the bilinear GF \eqref{mainthm}   of  Theorem \ref{mainthm-thm}. We refer the reader to \cite{chenwang} for the full proof. Furthermore, by comparing the coefficients of $t^n$, we also found

\begin{dl}[Cf. {\rm \cite[(3.2)]{chenwang}}]\label{hanid-II-dl}  Let $\Phi_{n}^{(\alpha)}(x, y | q) $ be given by \eqref{hanpoly}. Then for any integer $m\geq 0$, it holds
\begin{align}
\Phi_{n+m}^{(\alpha)}(x, y | q) =\sum_{k=0}^n \genfrac{[}{]}{0pt}{}{n}{k}_q
(\alpha x/y;q)_kx^{n-k}y^k\Phi_{m}^{(\alpha)}(x q^k,y| q).\label{hanid-II}
\end{align}
\end{dl}

Later as we will see,  this recurrence relation  sheds a new light on Rogers's symmetric expansion formula (cf. \cite[Thm. 1.1]{bowman-1}) and the $q$-Mehler formula \eqref{mainthm}. Actually,  it is this recurrence relation that we not only  show Rogers' symmetric expansion formula but also extend the $q$-Mehler formula \eqref{mainthm}  to a general version Theorem \ref{A-S-Cth1.1-genal} below, which is one of the main results of the present paper.

{\bf Notation:} For convenience of our discussions,  we will re-introduce the following  \textit{multivariate} Al-Salam-Carlitz polynomials, which are first posed in our paper \cite{chenwangma} but   also appeared in \cite[Sec. 1.1]{bowman-1} in a different form.
\begin{dy}[Cf. {\rm \cite[Def. 1.4]{chenwangma}}]\label{multi-hans}For $A:=(a_1,\cdots,a_m),  B:=(b_1,\cdots,b_m)\in\mathbb{C}^m$, $\mathbb{C}^m$ is the $m$-tuple vector space over the complex field $\mathbb{C}$,  we define the multivariate Al-Salam-Carlitz polynomials $\Phi_n(A;B)$ by
\begin{align}
\sum_{n=0}^{\infty}\Phi_n(A;B)\frac{t^n}{\poq{q}{n}}=\frac{\poq{At}{\infty}}{\poq{Bt}{\infty}},\label{geneal-hahn}
\end{align}
where, for clarity,  we write
$\poq{At}{n}$ for $\poq{a_1t,a_2t,\cdots,a_mt}{n}$ for integer $n\geq 0$ and suppress the  dimension $m$ and $q$. In particular, when $A=0$，i.e., $a_i=0$ for $1\leq i\leq m$, then we write $\Phi_n(B)$ for $\Phi_n(0;B)$   for simplicity.
\end{dy}
It is easily seen that
 $\Phi_n(A;B)$ are homogeneous, i.e., for any complex number $x$, it holds
 \begin{align}
\Phi_n(Ax;Bx)=x^n\Phi_n(A;B)\label{homo}
 \end{align}
 and can be given explicitly by the $q$-binomial theorem \cite[(II.3)]{10} in the form
\begin{align}\label{multi-ASC-ex}
    \Phi_n(A;B)=\sum_{i_1+i_2+\cdots+i_m=n}\genfrac{[}{]}{0pt}{}{n}{i_1, i_2,\ldots, i_m}_q\left[b_1, a_1\right]_{i_1} \cdots\left[b_m, a_m\right]_{i_m},
\end{align}
where the Cauchy polynomials $[x, y]_n$ are defined by $$\left[x, y\right]_{n}:=(x-y)(x-yq)\cdots(x-yq^{n-1})，$$
and the notation $\genfrac{[}{]}{0pt}{}{n}{i_1,i_2, \ldots, i_m}_q$ stands for the $q$-multinomial coefficient given by $$\frac{\poq{q}{n}}{\poq{q}{i_1}\poq{q}{i_2}\cdots\poq{q}{i_m}}.$$

For our convenience, we will write $\Phi_n(\alpha x;x,y)$ for $ \Phi_{n}^{(\alpha)}(x, y| q)$ above, since
\begin{align*}
\sum_{n=0}^{\infty}\Phi_{n}^{(\alpha)}(x, y| q) \frac{t^{n}}{(q ; q)_n}=\frac{(\alpha x t; q)_{\infty}}{(x t, yt; q)_{\infty}}.
\end{align*}

The present paper is organized  as follows. In the next section, along the same lines as the proof of Theorem \ref{temper} in our paper \cite{chenwang}, we will establish a general bilinear GF of the multivariate Al-Salam-Carlitz polynomials $\Phi_n(A;B)$.  Some applications to $q$-series closely related with $\Phi_n\left(A;B\right)$ are also investigated in Section \ref{sec3}. Two new symmetric expansion
formulas for ${}_2\phi_1$ and ${}_3\phi_2$ series being different from Rogers' results  are presented.

\section{Main theorems and proofs}\label{secfurther}
As planned, in this section we will establish a general bilinear GF of the multivariate Al-Salam-Carlitz polynomials $\Phi_n(A;B)$. Before this, we begin with a proof of Rogers' symmetric expansion, from which our main argument emerges.
\subsection{The proof of Rogers' symmetric expansion}
In his paper \cite[(13)]{rogers}, L. J. Rogers established the following symmetric expansion formula via the $q$-derivative operator method.
\begin{dl}[Cf. {\rm \cite[Thm. 1.1]{bowman-1}}]\label{fenmu-3-2}Let $\Phi_n\left(A;B\right)$ be given as above. Then it holds
\begin{align}&\sum_{n=0}^\infty \frac{\Phi_n\left(\lambda, \lambda_1\right)\Phi_n\left(\lambda_2, \lambda_3, \lambda_4\right) }{(q;q)_n}\label{fenmu-3-2-equ}\\
&=\frac{\left(\lambda \lambda_1 \lambda_2 \lambda_3;q\right)_{\infty}}{\left(\lambda \lambda_2,\lambda \lambda_3,\lambda_1 \lambda_2,\lambda_1 \lambda_3,\lambda_1 \lambda_4;q\right)_{\infty}}{ }_2 \phi_1\left[\begin{array}{c}
\lambda_1 \lambda_2, \lambda_1 \lambda_3 \\
\lambda \lambda_1 \lambda_2 \lambda_3
\end{array} ; q, \lambda \lambda_4\right] \nonumber.
\end{align}
\end{dl}
In the following, we will use Theorem \ref{hanid-II-dl} and merely by comparing coefficients, instead of the $q$-derivative operator,
to show Theorem \ref{fenmu-3-2}.

\pf It is clear that \eqref{fenmu-3-2-equ} can be reformulated as
    \begin{align*}
        \frac{1}{\poq{\lambda\lambda_2,\lambda\lambda_3}{\infty}}\sum_{n=0}^\infty \frac{\poq{\lambda \lambda_1 \lambda_2 \lambda_3q^n}{\infty}}{\poq{\lambda_1 \lambda_2q^n,\lambda_1 \lambda_3q^n}{\infty}\poq{q}{n}}\frac{(\lambda\lambda_4)^n}{\poq{\lambda_1\lambda_4}{\infty}}\\
        =\sum_{n=0}^\infty \Phi_n(\lambda,\lambda_1)\sum_{i+j+k=n}\frac{\lambda_2^i\lambda_3^j\lambda_4^k}{\poq{q}{i}\poq{q}{j}\poq{q}{k}}.
    \end{align*}
 In this form, we compare the coefficients of $\lambda_4^m$ on both sides, obtaining
    \begin{align*}
        \frac{1}{\poq{\lambda\lambda_2,\lambda\lambda_3}{\infty}}\sum_{n=0}^\infty \frac{\poq{\lambda \lambda_1 \lambda_2 \lambda_3q^n}{\infty}}{\poq{\lambda_1 \lambda_2q^n,\lambda_1 \lambda_3q^n}{\infty}}\frac{\lambda^n\lambda_1^{m-n}}{\poq{q}{n}\poq{q}{m-n}}\\
        =\sum_{n=0}^\infty\Phi_n(\lambda,\lambda_1)\sum_{i+j=n-m}\frac{\lambda_2^i\lambda_3^j}{\poq{q}{i}\poq{q}{j}\poq{q}{m}}.
    \end{align*}
  Now, by virtue of the $q$-Mehler formula \eqref{mainthm}  with $\alpha=\beta=0, t=1$, we find
    \begin{align*}
        \frac{\poq{\lambda \lambda_1 \lambda_2 \lambda_3q^n}{\infty}}{\poq{\lambda_1 \lambda_2q^n,\lambda_1 \lambda_3q^n,\lambda\lambda_2,\lambda\lambda_3}{\infty}}=\sum_{k=0}^\infty \frac{\Phi_k(\lambda_2,\lambda_3)
        \Phi_k(\lambda,\lambda_1q^n)}{\poq{q}{k}},
    \end{align*}
 which reduces the preceding identity to
    \begin{align*}
        \sum_{n=0}^\infty\frac{\lambda^n\lambda_1^{m-n}}{\poq{q}{n}\poq{q}{m-n}}\sum_{k=0}^\infty \frac{\Phi_k(\lambda_2,\lambda_3)\Phi_k(\lambda,\lambda_1q^n)}{\poq{q}{k}}
        =\sum_{n=0}^\infty\frac{\Phi_{n+m}(\lambda,\lambda_1)\Phi_n(\lambda_2,\lambda_3)}{\poq{q}{m}\poq{q}{n}}.
    \end{align*}
 Again, by comparing the coefficients of $\Phi_s(\lambda_2,\lambda_3)/\poq{q}{s}$ on both sides, we archive
    \begin{align*}
        \sum_{n=0}^m\genfrac{[}{]}{0pt}{}{m}{n}_q\Phi_s(\lambda,\lambda_1q^n)\lambda^n\lambda_1^{m-n}=\Phi_{s+m}(\lambda,\lambda_1).
    \end{align*}
 It is just the special case  $\alpha=0$  of  \eqref{hanid-II}. Thus Theorem \ref{fenmu-3-2} is proved.\qed

\subsection{A novel bilinear GF of $\Phi_n(A;B)$ and its variants}
Inspired by the argument of Theorem  \ref{fenmu-3-2}, we now proceed to set up a novel bilinear GF of $\Phi_n(A;B)$. To that end, we first need a preliminary which is key to our discussions.
\begin{yl}\label{added-yl-2}For $A:=(a_1,a_2,\cdots,a_m),  B:=(b_1,b_2,\cdots,b_m)\in \mathbb{C}^m$, and any integer $k\geq 0$, it holds
    \begin{align}
        \sum_{n=0}^\infty\Phi_{n+k}(A;B)\frac{x^{n}}{\poq{q}{n}}=\frac{\poq{Ax}{\infty}}{\poq{Bx}{\infty}}x^{-k}{}_{m+1}\phi_m\left[\begin{matrix}q^{-k},&Bx  \\ &Ax \end{matrix};q,q \right].\label{added-lemma}
    \end{align}
\end{yl}
\pf It is easy to check by definition (\ref{geneal-hahn}) that
\begin{align*}
    \sum_{n=0}^\infty\Phi_{n}(A;B)\frac{x^{n}}{\poq{q}{n}}=\frac{\poq{Ax}{\infty}}{\poq{Bx}{\infty}},
\end{align*}
thus, using the usual $q$-derivative operator $\mathcal{D}^k_{q,x}\, (k\geq 0)$
(cf. \cite[Exer. 1.12]{10}), it is easily found that
\begin{align*}
    \sum_{n=0}^\infty\Phi_{n}(A;B)\frac{\mathcal{D}^k_{q,x}(x^{n})}{\poq{q}{n}}=\mathcal{D}_{q,x}^k\bigg(\frac{\poq{Ax}{\infty}}{\poq{Bx}{\infty}}\bigg).
\end{align*}
It gives
\begin{align*}
    \sum_{n=0}^\infty\Phi_{n+k}(A;B)\frac{x^{n}}{\poq{q}{n}}&=x^{-k}\sum_{i=0}^k\genfrac{[}{]}{0pt}{}{k}{i}_q(-1)^i q^{\binom{i}{2}+i(1-k)}\frac{\poq{Axq^i}{\infty}}{\poq{Bxq^i}{\infty}}\\
    &=\frac{\poq{Ax}{\infty}}{\poq{Bx}{\infty}}x^{-k}\sum_{i=0}^k\frac{\poq{q^{-k},Bx}{i}}{\poq{q,Ax}{i}}q^i.
\end{align*}
That is we wanted.
\qed

Now we are in a good position to show a novel bilinear GF of $\Phi_n(A;B)$. As we will see later, it is indeed a generalization of
 Theorem \ref{mainthm-thm}.
\begin{dl}\label{A-S-Cth1.1-genal}Let $A,B\in \mathbb{C}^m$ and $\Phi_n(A;B)$ be the same as Definition \ref{multi-hans}. Then we have
   \begin{align} \sum_{n=0}^\infty&\Phi_{n}(\alpha ;x, y)\Phi_{n}(\beta, A;\lambda_2,\lambda_3,B)\frac{t^n}{\poq{q}{n}}\nonumber
   \\
   & =\frac{\poq{\alpha \lambda_3t,\beta yt,Axt}{\infty}}{\poq{x \lambda_3t, y\lambda_2t, y \lambda_3t,Bxt}{\infty}}\sum_{n=0}^\infty\frac{\poq{\alpha /y,x \lambda_3t}{n}}{\poq{q,\alpha \lambda_3 t}{n}}(y/x)^n\nonumber\\
   & \qquad \times{}_{m+1}\phi_m\left[\begin{matrix}q^{-n},Bxt  \\ Axt\end{matrix};q,q \right]{}_3\phi_2\left[\begin{matrix}
    \alpha/x,\beta/\lambda_2,y\lambda_3t\\
    \alpha \lambda_3tq^n,\beta yt
\end{matrix};q,x\lambda_2tq^n\right].\label{yyy-yyy-1}
    \end{align}
\end{dl}
\pf Let's start with \eqref{hanid-II}. We first reformulate it as form
\begin{align*}
    \frac{\Phi_{n+m}(\alpha x; x, y)}{\poq{q}{n}\poq{q}{m}}=\sum_{k=0}^n\frac{\poq{\alpha x/y}{k}}{\poq{q}{k}\poq{q}{n-k}}x^{n-k}y^k\frac{\Phi_{m}(\alpha xq^k;xq^k,y)}{\poq{q}{m}}.
\end{align*}
Next, multiply both sides by $\Phi_m(\beta \lambda_2;\lambda_2,\lambda_3)t^{m}$ and sum on $m$ from 0 to $+\infty$. The result is
\begin{align*}
    \sum_{m=0}^\infty&\frac{\Phi_{n+m}(\alpha x;x, y)\Phi_m(\beta \lambda_2;\lambda_2,\lambda_3)}{\poq{q}{n}\poq{q}{m}}t^{m}\\&=\sum_{k=0}^n\frac{\poq{\alpha x/y}{k}}{\poq{q}{k}\poq{q}{n-k}}x^{n-k}y^k\bigg(\sum_{m=0}^\infty\Phi_{m}(\alpha xq^k;xq^k,y)\Phi_m(\beta \lambda_2;\lambda_2,\lambda_3)\frac{t^{m}}{\poq{q}{m}}\bigg).
\end{align*}
Applying  the $q$-Mehler formula \eqref{mainthm} to the sum in brackets, we  have
\begin{align}\sum_{m=n}^\infty&\frac{\Phi_{m}(\alpha x;x, y)\Phi_{m-n}(\beta \lambda_2;\lambda_2,\lambda_3)}{\poq{q}{n}\poq{q}{m-n}}t^{m-n}=\sum_{k=0}^n\frac{\poq{\alpha x/y}{k}}{\poq{q}{k}\poq{q}{n-k}}x^{n-k}y^k\nonumber\\&\times\frac{(\alpha x\lambda_3tq^k, \beta y\lambda_2t ; q)_{\infty}}{(x \lambda_3tq^k, y \lambda_3t, y\lambda_2t ; q)_{\infty}}{}_3 \phi_2\left[\begin{array}{ccc}
\alpha, &\beta, &y\lambda_3t  \\
&\alpha x\lambda_3tq^k, &\beta y\lambda_2t
\end{array}; q, x\lambda_2tq^k\right].\label{11added}
\end{align}
Expand ${}_3\phi_2$ series in light of the definition and then interchange the order of the sums. We see that the right-hand side (in short, RHS) transforms into
\begin{align*}
\text{RHS of \eqref{11added}}
&=\frac{\poq{\alpha x\lambda_3t,\beta y\lambda_2t}{\infty}}{\poq{x \lambda_3t, y \lambda_3t, y\lambda_2t}{\infty}}\sum_{i=0}^\infty\frac{\poq{\alpha,\beta,y\lambda_3t}{i}}{\poq{q,\alpha x\lambda_3t,\beta y\lambda_2t}{i}}(x\lambda_2t)^i\nonumber\\
&\qquad\times\sum_{k=0}^n\frac{\poq{\alpha x/y,x \lambda_3t}{k}}{\poq{q,\alpha x\lambda_3tq^i}{k}}(yq^i)^k\frac{x^{n-k}}{\poq{q}{n-k}}.
\end{align*}
 Substituting this back to \eqref{11added} and then multiplying both sides with
 $\Phi_n(A;B)t^{n}$,  we are able to
  sum on $n$ from 0 to $+\infty$ and arrive at
\begin{align}
    \sum_{n=0}^\infty\sum_{m=n}^\infty&\frac{\Phi_{m}(\alpha x;x, y)\Phi_{m-n}(\beta \lambda_2;\lambda_2,\lambda_3)}{\poq{q}{n}\poq{q}{m-n}}\Phi_n(A;B)t^{m}\nonumber\\
    &=\frac{\poq{\alpha x\lambda_3t,\beta y\lambda_2t}{\infty}}{\poq{x \lambda_3t, y \lambda_3t, y\lambda_2t}{\infty}}\sum_{i=0}^\infty\frac{\poq{\alpha,\beta,y\lambda_3t}{i}}{\poq{q,\alpha x\lambda_3t,\beta y\lambda_2t}{i}}(x\lambda_2t)^i\nonumber\\
    &\quad \times\sum_{n=0}^\infty\sum_{k=0}^n\frac{\poq{\alpha x/y,x \lambda_3t}{k}}{\poq{q,\alpha x\lambda_3tq^i}{k}}(yq^i)^k\frac{x^{n-k}}{\poq{q}{n-k}}\Phi_n(A;B)t^{n}.\label{Phi-2-2-33}
    \end{align}
   Interchanging the order of the sums over $m$ and $n$ on the left-hand side (in short, LHS) of the above identity, and then evaluating by the relation (from the definition (\ref{geneal-hahn}))
\begin{align*}
    \sum_{n=0}^m\frac{\Phi_{m-n}(\beta \lambda_2; \lambda_2,\lambda_3)}{\poq{q}{m-n}}\frac{\Phi_n(A;B)}{\poq{q}{n}}=\frac{\Phi_m(\beta \lambda_2,A;\lambda_2,\lambda_3,B)}{\poq{q}{m}},
\end{align*}
we arrive at
\begin{align*}
    \text{LHS of \eqref{Phi-2-2-33}}
        &=\sum_{m=0}^\infty\Phi_{m}(\alpha x;x, y)\Phi_{m}(\beta \lambda_2, A;\lambda_2,\lambda_3,B)\frac{t^{m}}{\poq{q}{m}}.
\end{align*}
Summing up, we have
\begin{align*}
\sum_{m=0}^\infty&\Phi_{m}(\alpha x;x, y)\Phi_{m}(\beta \lambda_2, A;\lambda_2,\lambda_3,B)\frac{t^{m}}{\poq{q}{m}}\\
&=\frac{\poq{\alpha x\lambda_3t,\beta y\lambda_2t}{\infty}}{\poq{x \lambda_3t, y \lambda_3t, y\lambda_2t}{\infty}}\sum_{i=0}^\infty\frac{\poq{\alpha,\beta,y\lambda_3t}{i}}{\poq{q,\alpha x\lambda_3t,\beta y\lambda_2t}{i}}(x\lambda_2t)^i\\
    &\times\sum_{k=0}^\infty\frac{\poq{\alpha x/y,x \lambda_3t}{k}}{\poq{q,\alpha x\lambda_3tq^i}{k}}(ytq^i)^k\sum_{n=0}^\infty\Phi_{n+k}(A;B)\frac{(xt)^{n}}{\poq{q}{n}}.
    \end{align*}
A last replacement
$(\alpha,\beta)\to (\alpha/x, \beta/\lambda_2)$ and restated in terms of \eqref{added-lemma}
leads us to \eqref{yyy-yyy-1}. Thus the theorem is proved.
\qed

It is obvious that when $A=B$, \eqref{yyy-yyy-1} directly reduces to the $q$-Mehler formula \eqref{mainthm} under the simultaneous replacements $(\alpha,\beta)\to(\alpha x, \beta \lambda_2 )$. Apart from this,  Theorem \ref{A-S-Cth1.1-genal}  covers the following useful results when $t=1,A=(\gamma),B=(\lambda_1)$.
\begin{tl}\label{A-S-Cth1.1-together}
    \begin{align}
& \sum_{n=0}^{\infty} \frac{\Phi_n(\alpha ; x, y) \Phi_n\left(\beta, \gamma ; \lambda_1, \lambda_2, \lambda_3\right)}{(q ; q)_n} =\frac{\left(\alpha \lambda_3, \beta y, \gamma x ; q\right)_{\infty}}{\left(x \lambda_1, x \lambda_3, y \lambda_2, y \lambda_3 ; q\right)_{\infty}}\label{yyy-yyy-gel}\\
& \qquad\times \sum_{n=0}^{\infty} \frac{\left(\alpha / x, \beta / \lambda_2, y \lambda_3 ; q\right)_n}{\left(q, \alpha \lambda_3, \beta y ; q\right)_n}(x \lambda_2)^n\,{}_3 \phi_2\left[\begin{array}{c}
\alpha / y, x \lambda_3, \gamma / \lambda_1 \\
\alpha \lambda_3 q^n, \gamma x
\end{array} ; q, y \lambda_1 q^n\right].\nonumber
    \end{align}
 Particularly, when $\gamma \to 0$,  we have
\begin{align}        &\sum_{n=0}^\infty\frac{\Phi_{n}(\alpha ;x, y)\Phi_{n}(\beta ;\lambda_1, \lambda_2,\lambda_3)}{\poq{q}{n}}=\frac{\poq{\alpha \lambda_3,\beta y}{\infty}}{\poq{x \lambda_1,x\lambda_3, y \lambda_2, y\lambda_3}{\infty}}\label{yyy-yyy}\\&\qquad\times\sum_{n=0}^\infty\frac{\poq{\alpha/x,\beta/\lambda_2,y\lambda_3}{n}}{\poq{q,\alpha \lambda_3,\beta y}{n}}(x\lambda_2)^n\,{}_2\phi_1\left[\begin{matrix}
        \alpha/y,x \lambda_3\\\alpha \lambda_3q^n   \end{matrix};q,y\lambda_1q^n\right].\nonumber
    \end{align}
\end{tl}

The special case $\alpha=x$ and $t=1$ of \eqref{yyy-yyy-1} deserves a separate statement, from which  we may derive a generalization of the $q$-Gauss sum \cite[(II.8)]{10}.
    \begin{tl}\label{A-S-Cth1.1-genal-one}For any $A,B\in \mathbb{C}^m$ and $x,y$ with $|y/x|<1$, we have
  \begin{align}
       \sum_{n=0}^\infty\frac{\poq{x /y}{n}}{\poq{q}{n}}(y/x)^n{}_{m+1}\phi_m\left[\begin{matrix}q^{-n},&Bx  \\ &Ax\end{matrix};q,q \right]
    =\frac{\poq{Ay,Bx}{\infty}}{\poq{By,Ax}{\infty}}. \label{interest}
    \end{align}
\end{tl}
\pf We show \eqref{interest} in two ways.

{\bf\sl (The first proof)}. Let $\alpha=x, t=1$ in \eqref{yyy-yyy-1}. We have
 \begin{align*}
    \sum_{n=0}^\infty\frac{\Phi_{n}(\beta, A ;\lambda_2,\lambda_3, B)}{\poq{q}{n}}y^n
    &=\frac{\poq{\beta y,Ax}{\infty}}{\poq{y\lambda_2, y \lambda_3,Bx}{\infty}}\\
    &\times\sum_{k=0}^\infty\frac{\poq{x /y}{k}}{\poq{q}{k}}(y/x)^k{}_{m+1}\phi_m\left[\begin{matrix}q^{-k},&Bx  \\ &Ax\end{matrix};q,q \right].
    \end{align*}
  On the other hand,   by the definition of $\Phi_{n}$, it also holds
    \begin{align*}
    \sum_{n=0}^\infty\frac{\Phi_{n}(\beta, A;\lambda_2,\lambda_3, B)}{\poq{q}{n}}y^n
    =\frac{\poq{\beta y, Ay}{\infty}}{\poq{\lambda_2 y,\lambda_3 y, By}{\infty}}.
    \end{align*}
  Collecting both together, we  get \eqref{interest}.

 ({\bf\sl The second proof}). To show  \eqref{interest}, let $y=tx$. Then we only need to check that for $|t|<1$, it holds
\begin{align*}
       \sum_{k=0}^\infty\frac{\poq{1/t}{k}}{\poq{q}{k}}t^k{}_{m+1}\phi_m\left[\begin{matrix}q^{-k},&Bx  \\ &Ax\end{matrix};q,q \right]
    =\frac{\poq{Atx,Bx}{\infty}}{\poq{Btx,Ax}{\infty}}.
    \end{align*}
    By Ismail's argument (see \cite{schlosser}) or analytic continuation, it suffices to check it is valid for $t=q^{N},N\geq 0$,
 \begin{align*}
       \sum_{k=0}^N\genfrac{[}{]}{0pt}{}{N}{k}_q\tau(k){}_{m+1}\phi_m\left[\begin{matrix}q^{-k},&Bx  \\ &Ax\end{matrix};q,q \right]
    =\frac{\poq{Bx}{N}}{\poq{Ax}{N}}.
    \end{align*}
    By inverting, it becomes
 \begin{align*}
       \sum_{k=0}^N\genfrac{[}{]}{0pt}{}{N}{k}_q\tau(N-k)\frac{\poq{Bx}{k}}{\poq{Ax}{k}}
    =\tau(N){}_{m+1}\phi_m\left[\begin{matrix}q^{-N},&Bx  \\ &Ax\end{matrix};q,q \right].
    \end{align*}
 That is
\begin{align*}
       \sum_{k=0}^N\frac{\poq{q^{-N},Bx}{k}}{\poq{q,Ax}{k}}q^k
    ={}_{m+1}\phi_m\left[\begin{matrix}q^{-N},&Bx  \\ &Ax\end{matrix};q,q \right].
\end{align*}
    It is self-evident. Thus \eqref{interest} is confirmed.
\qed

As mentioned before, \eqref{interest} indeed contains the
 $q$-Gauss sum \cite[(II.8)]{10} as the special case $A=(c),B=(d)$, namely
 \begin{align*}
      {}_2\phi_1\left[\begin{matrix}
          x/y,c/d\\
          cx
      \end{matrix}q;dy\right]
    =\frac{\poq{cy,dx}{\infty}}{\poq{dy,cx}{\infty}}.
    \end{align*}

We close this section by another interesting case  $\beta=\lambda_2$ and $t=1$ of \eqref{yyy-yyy-1} in Theorem \ref{A-S-Cth1.1-genal}, which seems more general than \eqref{interest}.
\begin{tl} For $x,y$ with $|y/x|<1$, we have
\begin{align}        &\sum_{n=0}^\infty\frac{\Phi_{n}(\alpha ;x, y)\Phi_{n}(A;\lambda_3, B)}{\poq{q}{n}}\label{yyy-yyy-1-2}\\
        &=\frac{\poq{\alpha \lambda_3,Ax}{\infty}}{\poq{x \lambda_3, y \lambda_3,Bx}{\infty}}\sum_{k=0}^\infty\frac{\poq{\alpha /y,x \lambda_3}{k}}{\poq{q,\alpha \lambda_3}{k}}(y/x)^k\,{}_{m+1}\phi_m\left[\begin{matrix}q^{-k},&Bx  \\ &Ax\end{matrix};q,q \right]. \nonumber
    \end{align}
\end{tl}

\section{Applications}\label{sec3}
This part is devoted to some concrete GFs derived  from  Theorem \ref{A-S-Cth1.1-genal} and Corollary \ref{A-S-Cth1.1-together} under proper series transformations.
\subsection{Results with the Rogers-Hall transformation involved}
By virtue of Rogers-Hall's second transformation \cite[(III.10)]{10}, it is easy to
check
\begin{align*}
   &{}_3\phi_2\left[\begin{matrix}
        \alpha/y,\gamma/\lambda_1,x \lambda_3\\\alpha \lambda_3q^n, x\gamma   \end{matrix};q,y\lambda_1q^n\right]\\
       = &\frac{\poq{\gamma/\lambda_1, x y\lambda_1 \lambda_3  q^n, \alpha \lambda_1 q^n }{\infty}}{\poq{ x\gamma, \alpha \lambda_3 q^n, y\lambda_1 q^n }{\infty}}{}_3\phi_2\left[\begin{matrix}
       x\lambda_1 , \alpha \lambda_1 \lambda_3 q^n/\gamma, y\lambda_1 q^n\\ x y\lambda_1 \lambda_3 q^n, \alpha \lambda_1 q^n \end{matrix};q, \gamma/\lambda_1\right].
\end{align*}
Upon applying  this transformation to \eqref{yyy-yyy-gel}, we easily establish
\begin{dl}\label{A-S-Cth1.1-22}Let $\Phi_n(A;B)$ be given as above. Then    \begin{align}        &\sum_{n=0}^\infty\frac{\Phi_{n}(\alpha ;x, y)\Phi_{n}(\beta, \gamma ;\lambda_1, \lambda_2,\lambda_3)}{\poq{q}{n}}=\frac{\poq{\alpha \lambda_1,\beta y,\gamma/\lambda_1,x y\lambda_1 \lambda_3}{\infty}}{\poq{x \lambda_1,x\lambda_3,y \lambda_1, y \lambda_2, y\lambda_3}{\infty}}\label{yyy-yyy-ooo}\\
&\qquad\times\sum_{n=0}^\infty\frac{\poq{\alpha/x,\beta/\lambda_2, y\lambda_1,y\lambda_3}{n}}{\poq{q,\alpha \lambda_1,\beta y,x y\lambda_1 \lambda_3}{n}}(x\lambda_2)^n\,{}_3\phi_2\left[\begin{matrix}
       x\lambda_1 , \alpha \lambda_1 \lambda_3 q^n/\gamma, y\lambda_1 q^n \\  x y\lambda_1 \lambda_3 q^n, \alpha \lambda_1 q^n \end{matrix};q, \gamma/\lambda_1\right]\nonumber.
    \end{align}
\end{dl}

    Putting $\gamma\to 0$, then we obtain
    \begin{tl}\label{A-S-Cth1.1-22-33}
    Let $\Phi_n(A;B)$ be given as above. Then
    \begin{align}        &\sum_{n=0}^\infty\frac{\Phi_{n}(\alpha ;x, y)\Phi_{n}(\beta;\lambda_1, \lambda_2,\lambda_3)}{\poq{q}{n}}=\frac{\poq{\alpha \lambda_1,\beta y,x y\lambda_1 \lambda_3}{\infty}}{\poq{x \lambda_1,x\lambda_3,y \lambda_1, y \lambda_2, y\lambda_3}{\infty}}\label{yyy-yyy-ooo-000}\\
&\qquad\times \sum_{n=0}^\infty\frac{\poq{\alpha/x,\beta/\lambda_2,y \lambda_1,y\lambda_3}{n}}{\poq{q,\alpha \lambda_1,\beta y,x y\lambda_1 \lambda_3}{n}}(x\lambda_2)^n\,{}_2\phi_2\left[\begin{matrix}
       x\lambda_1 , y\lambda_1 q^n\\  x y\lambda_1 \lambda_3 q^n, \alpha \lambda_1 q^n \end{matrix};q, \alpha\lambda_3 q^n\right]\nonumber.
    \end{align}
\end{tl}

The following is another generalization of Bowman's result  (cf. \cite[Thm. 3.11]{bowman-1}).
\begin{tl} Let $\Phi_n(A;B)$ be given as above. Then  we have
   \begin{align}    &\sum_{n=0}^\infty\frac{\Phi_{n}(\alpha ;x, y)\Phi_{n}(\beta;\lambda_2,\lambda_3,\lambda_4,\lambda_5)}{\poq{q}{n}}\label{yyy-yyy-1-123}
    =\frac{\poq{\alpha \lambda_3,\beta y}{\infty}}{\poq{x \lambda_3, y\lambda_2, y \lambda_3,x\lambda_4,x\lambda_5}{\infty}}\\
   &\qquad  \times\sum_{n=0}^\infty\frac{\poq{\alpha /y,x \lambda_3}{n}}{\poq{q,\alpha \lambda_3}{n}}y^n\Phi_n(x\lambda_4\lambda_5;\lambda_4,\lambda_5){}_3\phi_2\left[\begin{matrix}
    \alpha/x,\beta/\lambda_2,y\lambda_3\\
    \alpha \lambda_3q^n,\beta y
\end{matrix};q,x\lambda_2q^n\right]\nonumber.
\end{align}
\end{tl}
\pf  It suffices  to set $t=1, A=(0,0),B=(\lambda_4,\lambda_5)\in \mathbb{C}^2$ in \eqref{yyy-yyy-1}. Then we have
 \begin{align}    \sum_{n=0}^\infty\frac{\Phi_{n}(\alpha ;x, y)\Phi_{n}(\beta;\lambda_2,\lambda_3,\lambda_4,\lambda_5)}{\poq{q}{n}}\label{yyy-yyy-1-1234}
    =\frac{\poq{\alpha \lambda_3,\beta y}{\infty}}{\poq{x \lambda_3, y\lambda_2, y \lambda_3,x\lambda_4,x\lambda_5}{\infty}}\\
    \times\sum_{n=0}^\infty\frac{\poq{\alpha /y,x \lambda_3}{n}}{\poq{q,\alpha \lambda_3}{n}}(y/x)^n{}_{3}\phi_2\left[\begin{matrix}q^{-n},x\lambda_4,x\lambda_5  \\ 0,0\end{matrix};q,q \right]\nonumber\\
    \times{}_3\phi_2\left[\begin{matrix}
    \alpha/x,\beta/\lambda_2,y\lambda_3\\
    \alpha \lambda_3q^n,\beta y
\end{matrix};q,x\lambda_2q^n\right]\nonumber.
\end{align}
Now, by applying the transformation \cite[(III.12)]{10}, we get
\begin{align*}
    {}_{3}\phi_2\left[\begin{matrix}q^{-n},&x\lambda_4,&x\lambda_5  \\ &0,&0\end{matrix};q,q \right]&=\lim_{d,e\to 0}{}_{3}\phi_2\left[\begin{matrix}q^{-n},&x\lambda_4,&x\lambda_5  \\ &dx,&ex\end{matrix};q,q \right]\\
    &=\lim_{d,e\to 0}\frac{(e /\lambda_5 ; q)_n}{(e ; q)_n} \lambda_5^n~{ }_3 \phi_2\left[\begin{array}{l}
q^{-n}, x\lambda_5, d /\lambda_4 \\
d, \lambda_5 q^{1-n} / e
\end{array} ; q,\frac{x\lambda_4 q}{e}\right]\\
    &=(x\lambda_5)^n\sum_{k=0}^n\frac{\poq{q^{-n},x\lambda_5}{k}}{\poq{q}{k}}\frac{(\lambda_4q^{n}/\lambda_5)^k}{\tau(k)}\\
    &=x^n\sum_{k=0}^n\genfrac{[}{]}{0pt}{}{n}{k}_q\poq{x\lambda_5}{k}\lambda_4^k\lambda_5^{n-k}=x^n\Phi_n(x\lambda_4\lambda_5;\lambda_4,\lambda_5).
\end{align*}
Inserting this back to \eqref{yyy-yyy-1-1234}, we obtain  \eqref{yyy-yyy-1-123} at once.
\qed
\begin{remark} Indeed, when $\alpha=\beta=0$ in \eqref{yyy-yyy-1-123},  we recover   \cite[Thm. 3.11]{bowman-1} due to Bowman as below. Bowman obtained this result by the $q$-derivative operator along Rogers' idea.
 \begin{align*}
     \sum_{n=0}^\infty\frac{\Phi_{n}(x, y)\Phi_{n}(\lambda_2,\lambda_3,\lambda_4,\lambda_5)}{\poq{q}{n}}    &=\frac{\poq{xy\lambda_2\lambda_3}{}}{\poq{x\lambda_2,x \lambda_3, x\lambda_4,x\lambda_5,y\lambda_2, y \lambda_3}{\infty}}\\&\times\sum_{n=0}^\infty\frac{\poq{x\lambda_2,x \lambda_3}{n}}{\poq{q,xy\lambda_2\lambda_3}{n}}\Phi_n(x\lambda_4\lambda_5;\lambda_4,\lambda_5)y^n.
\end{align*}
Alternately, putting $\alpha=x$ in \eqref{yyy-yyy-1-123}, then we may derive
\begin{align*} \sum_{n=0}^\infty\Phi_n(x\lambda_4\lambda_5;\lambda_4,\lambda_5)\frac{\poq{x/y}{n}}{\poq{q}{n}}y^n=\frac{\poq{x\lambda_4,x\lambda_5}{\infty}}{\poq{y\lambda_4, y \lambda_5}{\infty}}\nonumber.
\end{align*}
\end{remark}

\subsection{Results with Heine's transformations used}
In what follows, we will focus on applications of both \eqref{yyy-yyy-gel} and \eqref{yyy-yyy} in Corollary \ref{A-S-Cth1.1-together}.

\subsubsection{Using Heine's second transformation}
From Heine's second transformation \cite[(III.2)]{10}, it follows
\begin{align}
    {}_2\phi_1\left[\begin{matrix}
        \alpha/y,x \lambda_3\\\alpha \lambda_3q^n    \end{matrix};q,y\lambda_1q^n\right]
        =\frac{(\alpha q^n/x,xy\lambda_1\lambda_3q^n;q)_\infty}{(\alpha\lambda_3q^n,y\lambda_1 q^n;q)_\infty}{}_2\phi_1\left[\begin{matrix}
        x\lambda_1,x \lambda_3\\xy\lambda_1\lambda_3q^n    \end{matrix};q,\alpha q^n/x\right],\label{www-www}
\end{align}
which reduces \eqref{yyy-yyy} to
\begin{dl}
\begin{align}        \sum_{n=0}^\infty\frac{\Phi_{n}(\alpha x ;x, y)\Phi_{n}(\beta ;\lambda_1, \lambda_2,\lambda_3)}{\poq{q}{n}}\label{zzz}
          =\frac{\poq{\alpha,\beta y,xy\lambda_1\lambda_3}{\infty}}{\poq{x \lambda_1,x\lambda_3, y \lambda_1, y\lambda_2,y\lambda_3}{\infty}}\\
          \times\sum_{n=0}^\infty\frac{\poq{\beta/\lambda_2,y\lambda_1,y\lambda_3}{n}}{\poq{q,\beta y,xy\lambda_1\lambda_3}{n}}(x\lambda_2)^n\,{}_2\phi_1\left[\begin{matrix}
        x\lambda_1,x \lambda_3\\xy\lambda_1\lambda_3q^n    \end{matrix};q,\alpha q^n\right]\nonumber.    \end{align}
\end{dl}
\pf
 It only needs to check  directly, by substituting \eqref{www-www} into, that
\begin{align*}        \sum_{n=0}^\infty&\frac{\Phi_{n}(\alpha ;x, y)\Phi_{n}(\beta ;\lambda_1, \lambda_2,\lambda_3)}{\poq{q}{n}}\\
&=\frac{\poq{\alpha \lambda_3,\beta y}{\infty}}{\poq{x \lambda_1,x\lambda_3, y \lambda_2, y\lambda_3}{\infty}}\sum_{n=0}^\infty\frac{\poq{\alpha/x,\beta/\lambda_2,y\lambda_3}{n}}{\poq{q,\alpha \lambda_3,\beta y}{n}}(x\lambda_2)^n\\
&\qquad \times \frac{(\alpha q^n/x,xy\lambda_1\lambda_3q^n;q)_\infty}{(\alpha\lambda_3q^n,y\lambda_1 q^n;q)_\infty}{}_2\phi_1\left[\begin{matrix}
        x\lambda_1,x \lambda_3\\xy\lambda_1\lambda_3q^n    \end{matrix};q,\alpha q^n/x\right]\\
        &=\frac{\poq{\alpha/x,\beta y,xy\lambda_1\lambda_3}{\infty}}{\poq{x \lambda_1,x\lambda_3, y \lambda_2, y\lambda_3,y\lambda_1}{\infty}}\sum_{n=0}^\infty\frac{\poq{\beta/\lambda_2,y\lambda_1,y\lambda_3}{n}}{\poq{q,\beta y,xy\lambda_1\lambda_3}{n}}(x\lambda_2)^n\\
        &\qquad \times{}_2\phi_1\left[\begin{matrix}
        x\lambda_1,x \lambda_3\\xy\lambda_1\lambda_3q^n    \end{matrix};q,\alpha q^n/x\right]\nonumber.
    \end{align*}
  It is in agreement with \eqref{zzz} after replacing $\alpha$ with $\alpha x$.
\qed

Once taking $\alpha=0$ in \eqref{zzz}, we come up with the following
common extension of \cite[Thms. 1.1, 1.12, 2.1; Cor. 3.17]{bowman-1}. We leave all verifications to the interested reader.
\begin{xz}
\begin{align*}        &\sum_{n=0}^\infty\frac{\Phi_{n}(x, y)\Phi_{n}(\beta ;\lambda_1, \lambda_2,\lambda_3)}{\poq{q}{n}}
\\
          &=\frac{\poq{\beta y,xy\lambda_1\lambda_3}{\infty}}{\poq{x \lambda_1,x\lambda_3, y \lambda_1, y\lambda_2,y\lambda_3}{\infty}}{}_3\phi_2\left[\begin{matrix}
       \beta/\lambda_2,y\lambda_1,y\lambda_3\\\beta y,xy\lambda_1\lambda_3    \end{matrix};q,x\lambda_2\right]\nonumber.
       \end{align*}
\end{xz}
Instead, by taking $y=0$ in \eqref{zzz}, we obtain
\begin{xz}
 \begin{align*}
 \sum_{n=0}^\infty\Phi_{n}(\beta ;\lambda_1, \lambda_2,\lambda_3)\frac{\poq{\alpha}{n}}{\poq{q}{n}}x^n=\frac{\poq{\alpha,x\beta }{\infty}}{\poq{x \lambda_1,x \lambda_2,x\lambda_3}{\infty}}{}_3\phi_2\left[\begin{matrix}
         x\lambda_1,x \lambda_2, x \lambda_3\\\beta x,0   \end{matrix};q,\alpha\right].
         \end{align*}
\end{xz}
\pf Note that in this situation, it holds
\(\Phi_{n}(\alpha x;x, 0)=[x,x\alpha]_n=\poq{\alpha}{n}x^n,\)
thus gives rise to
\begin{align*}        &\sum_{n=0}^\infty\frac{[x,x\alpha]_n\Phi_{n}(\beta ;\lambda_1, \lambda_2,\lambda_3)}{\poq{q}{n}}\\
         &=\frac{\poq{\alpha}{\infty}}{\poq{x \lambda_1,x\lambda_3}{\infty}}\sum_{k=0}^\infty\frac{(x\lambda_1,x \lambda_3;q)_k}{\poq{q}{k}}  \alpha^k
         \sum_{n=0}^\infty \frac{\poq{\beta/\lambda_2}{n}}{\poq{q}{n}}(x\lambda_2q^k)^n \\
          &=\frac{\poq{\alpha}{\infty}}{\poq{x \lambda_1,x\lambda_3}{\infty}}\sum_{k=0}^\infty\frac{(x\lambda_1,x \lambda_3;q)_k}{\poq{q}{k}}  \alpha^k   \frac{\poq{\beta xq^k}{\infty}}{\poq{x\lambda_2q^k}{\infty}}\\
          &=\frac{\poq{\alpha,\beta x}{\infty}}{\poq{x \lambda_1,x \lambda_2,x\lambda_3}{\infty}}\sum_{k=0}^\infty\frac{(x\lambda_1,x \lambda_2, x \lambda_3;q)_k}{\poq{q,\beta x}{k}}  \alpha^k.
        \end{align*}
In this process, we have used the $q$-binomial theorem \cite[(II.3)]{10}. It is just the desired result.
\qed
\subsubsection{Using Heine's third transformation}
 From Heine's third transformation \cite[(III.3)]{10}, it follows
\begin{align}
    {}_2\phi_1\left[\begin{matrix}
        \alpha/y,x \lambda_3\\\alpha \lambda_3q^n    \end{matrix};q,y\lambda_1q^n\right]=\frac{(x\lambda_1;q)_\infty}{(y\lambda_1q^n;q)_\infty}{}_2\phi_1\left[\begin{matrix}
        \alpha q^n/x,y \lambda_3 q^n\\\alpha \lambda_3q^n    \end{matrix};q,x\lambda_1\right],\label{www-www-123}
\end{align}which helps us to reduce \eqref{yyy-yyy} to
\begin{dl}
\begin{align}        &\sum_{n=0}^\infty\frac{\Phi_{n}(\alpha ;x, y)\Phi_{n}(\beta ;\lambda_1, \lambda_2,\lambda_3)}{\poq{q}{n}}=\frac{\poq{\alpha \lambda_3,\beta y}{\infty}}{\poq{x \lambda_3, y \lambda_1, y\lambda_2,y\lambda_3}{\infty}}\label{zzz-new}\\
          &\qquad \times\sum_{n=0}^\infty\frac{\poq{\alpha/x,y\lambda_3}{n}}{\poq{q,\alpha \lambda_3}{n}}(x\lambda_1)^n\sum_{k=0}^n\genfrac{[}{]}{0pt}{}{n}{k}_q\frac{        (\beta/\lambda_2,y\lambda_1;q)_k}{(\beta y;q)_k }(\lambda_2/\lambda_1)^{k}\nonumber.
    \end{align}
\end{dl}
\pf
 Proceeding as above, we check in a straightforward way by  \eqref{www-www-123} that
\begin{align*}        \sum_{n=0}^\infty\frac{\Phi_{n}(\alpha ;x, y)\Phi_{n}(\beta ;\lambda_1, \lambda_2,\lambda_3)}{\poq{q}{n}}&\\
        =\frac{\poq{\alpha \lambda_3,\beta y}{\infty}}{\poq{x \lambda_3, y\lambda_1,y \lambda_2, y\lambda_3}{\infty}}\sum_{k=0}^\infty&\frac{\poq{\alpha/x,\beta/\lambda_2,y\lambda_1,y\lambda_3}{k}}{\poq{q,\alpha \lambda_3,\beta y}{k}}(x\lambda_2)^k\\&\times{}_2\phi_1\left[\begin{matrix}
        \alpha q^k/x,y \lambda_3 q^k\\\alpha \lambda_3q^k    \end{matrix};q,x\lambda_1\right]\nonumber\\
    =\frac{\poq{\alpha \lambda_3,\beta y}{\infty}}{\poq{x \lambda_3, y \lambda_1, y\lambda_2,y\lambda_3}{\infty}}\sum_{n=0}^\infty&\frac{\poq{\alpha/x,y\lambda_3}{n}}{\poq{q,\alpha \lambda_3}{n}}x^n\sum_{k=0}^n\genfrac{[}{]}{0pt}{}{n}{k}_q\frac{        (\beta/\lambda_2,y\lambda_1;q)_k}{(\beta y;q)_k }\lambda_2^k\lambda_1^{n-k}.
    \end{align*}
This completes the proof. \qed

Let $\lambda_2\to 0$ in \eqref{zzz-new}. Then we recover Theorem \ref{mainthm-thm} at once
\begin{lz}
\begin{align}        \sum_{n=0}^\infty\frac{\Phi_{n}(\alpha ;x, y)\Phi_{n}(\beta ; \lambda_1,\lambda_3)}{\poq{q}{n}}
          =\frac{\poq{\alpha \lambda_3,\beta y}{\infty}}{\poq{x \lambda_3,  y\lambda_1,y\lambda_3}{\infty}}{}_3\phi_2\left[\begin{matrix}
        \alpha/x,y\lambda_3,\beta/\lambda_1\\\alpha \lambda_3,\beta y   \end{matrix};q,x\lambda_1\right]\label{zzz-uuu}.
    \end{align}
\end{lz}

Let $\beta=\lambda_1$ in \eqref{zzz-new}. Then we obtain
\begin{lz}
\begin{align}        \sum_{n=0}^\infty\frac{\Phi_{n}(\alpha ;x, y)\Phi_{n}(\lambda_2, \lambda_3)}{\poq{q}{n}}=\frac{\poq{\alpha \lambda_3}{\infty}}{\poq{x \lambda_3, y \lambda_2, y\lambda_3}{\infty}}{}_2\phi_1\left[\begin{matrix}
        \alpha/x,y\lambda_3\\\alpha \lambda_3    \end{matrix};q,x\lambda_2\right].\label{yyy-yyy-yyy1}
    \end{align}
 \end{lz}
 \pf It follows from \cite[Exer. 1.3]{10} that
  \begin{align*}(a b ; q)_n=\sum_{k=0}^n
  \genfrac{[}{]}{0pt}{}{n}{k}_q b^k(a ; q)_k(b ; q)_{n-k}
  \end{align*}
  with $a=0$ and $b=\lambda_1/\lambda_2$, yielding
 \begin{align}        \sum_{k=0}^n\genfrac{[}{]}{0pt}{}{n}{k}_q    (\lambda_1/\lambda_2;q)_k\lambda_2^k\lambda_1^{n-k}=\lambda_2^n\nonumber.
    \end{align}
    This reduces \eqref{zzz-new} to \eqref{yyy-yyy-yyy1}.
 \qed

Using \eqref{yyy-yyy-yyy1}, we may find another expansion formula for ${}_2\phi_1$ series which is evidently different from Rogers' expansion.
    \begin{dl}
    \begin{align} {}_2\phi_1\left[\begin{matrix}
        a,b\\c    \end{matrix};q,d\right]= \frac{\poq{b,c/a,
        a b d/c}{\infty}}{\poq{c}{\infty}}\sum_{n=0}^\infty\frac{\Phi_{n}(a; 1, a b/c)\Phi_{n}(d, c /a)}{\poq{q}{n}}. \label{finalid}
    \end{align}
    \end{dl}
    \pf  In fact, by making the replacements in \eqref{yyy-yyy-yyy1}
    \[(\alpha,x,y,\lambda_3)\to\bigg(\frac{a d}{\lambda_2},\frac{d}{\lambda_2},\frac{a b d}{c \lambda_2},\frac{c\lambda_2}{a d}\bigg),\]
   we have
    \begin{align*}  {}_2\phi_1\left[\begin{matrix}
        a,b\\c    \end{matrix};q,d\right]&=\frac{\poq{b,c/a, a b d/c}{\infty}}{\poq{c}{\infty}}\nonumber\\
        &\times\sum_{n=0}^\infty\frac{\Phi_{n}(a d/\lambda_2; d/\lambda_2, a b d/(c \lambda_2))\Phi_{n}(\lambda_2, c \lambda_2/(a d))}{\poq{q}{n}}\nonumber.
    \end{align*}
   It is easy to verify by use of \eqref{homo}  that
   \begin{align*}
   \Phi_{n}(a d/\lambda_2; d/\lambda_2, a b d/(c \lambda_2))&=(d/\lambda_2)^{n}\Phi_{n}(a; 1, a b/c),\\
   \Phi_{n}(\lambda_2, c \lambda_2/(a d))&=(\lambda_2/d)^n\Phi_{n}(d, c /a).
    \end{align*}
           All these lead us to
     \begin{align*} {}_2\phi_1\left[\begin{matrix}
        a,b\\c    \end{matrix};q,d\right]
        = \frac{\poq{b,c/a, a b d/c}{\infty}}{\poq{c}{\infty}}\sum_{n=0}^\infty\frac{\Phi_{n}(a; 1, a b/c)\Phi_{n}(d, c /a)}{\poq{q}{n}}.
    \end{align*}
    The theorem is proved.
    \qed

    \begin{lz}
         If $d=0$ in \eqref{finalid}, then it becomes
    \begin{align*}
    \sum_{n=0}^\infty\frac{\Phi_{n}(a; 1,a b /c)}{\poq{q}{n}}\bigg(\frac{c}{a}\bigg)^n=\frac{\poq{c}{\infty}}{\poq{c/a,  b}{\infty}}.
    \end{align*}
    \end{lz}
\begin{lz}Let $b=c$ in \eqref{finalid}. Then  by the $q$-binomial theorem \cite[(II.3)]{10}, we have
    \begin{align*} \sum_{n=0}^\infty\frac{\Phi_{n}(d, b /a)}{\poq{q}{n}}=\frac{1}{\poq{b/a, d}{\infty}}.
    \end{align*}
    Note that $\Phi_{n}(a; 1, a)=1$.
    \end{lz}
\section{New symmetric expansions for ${}_2\phi_{1}$ and ${}_3\phi_{2}$ series}
Undoubtedly, Rogers' symmetric expansion formula \eqref{fenmu-3-2-equ} answered the problem of existence of Heine's transformations for ${}_2\phi_1$ series. See \cite[Sec. 1]{bowman-1} for more details. Thus, it is worthwhile to compare  our foregoing expansion formulas  with Rogers' symmetric expansion formula. For this purpose, we  make the simultaneous replacements
\[(a,b,c,d)\to (\lambda_1 \lambda_2, \lambda_1 \lambda_3,
\lambda \lambda_1 \lambda_2 \lambda_3, \lambda \lambda_4)\]
in the expansion  \eqref{finalid} and then reformulate it as
\begin{dl}\label{fenmu-3-2-222}Let $\Phi_n\left(A;B\right)$ be defined as above. Then it holds
\begin{align}&{ }_2 \phi_1\left[\begin{array}{c}
\lambda_1 \lambda_2, \lambda_1 \lambda_3 \\
\lambda \lambda_1 \lambda_2 \lambda_3
\end{array} ; q, \lambda \lambda_4\right]\label{fenmu-3-2-equequ}\\
&=\frac{\left(\lambda \lambda_3,\lambda_1 \lambda_3,\lambda_1 \lambda_4;q\right)_{\infty}} {\left(\lambda \lambda_1 \lambda_2 \lambda_3;q\right)_{\infty}}\sum_{n=0}^\infty \frac{\Phi_n\left(\lambda \lambda_1 \lambda_2;\lambda, \lambda_1\right)\Phi_n\left(\lambda_3,  \lambda_4\right) }{(q;q)_n}\nonumber.
\end{align}
\end{dl}
It is worth pointing out  that  \eqref{fenmu-3-2-equequ} is different from Rogers' symmetric expansion \eqref{fenmu-3-2-equ}. In the meantime, we are able to recover three Heine's  transformations which follow, respectively, from the invariance of the RHS of \eqref{fenmu-3-2-equequ} with respect to  the subgroup $$S_2\times S_2=<(\lambda_3,\lambda_4)\to (\lambda_4,\lambda_3), (\lambda,\lambda_1)\to (\lambda_1,\lambda)>,$$ not $S_2\times S_3$ usually known to the literature. Here, $S_n$ denotes the symmetric group on $n$ elements. Regarding this automorphism group we refer the reader  to \cite[pp. 1--2]{bowman-1}.  In addition, a combination \eqref{fenmu-3-2-equequ}  with \eqref{fenmu-3-2-equ} yields
\begin{tl}
\begin{align}
   \frac{1}{\left(\lambda \lambda_2,\lambda_1 \lambda_2;q\right)_{\infty}}
\sum_{n=0}^\infty \frac{\Phi_n\left(\lambda \lambda_1 \lambda_2;\lambda, \lambda_1\right)\Phi_n\left(\lambda_3,  \lambda_4\right) }{(q;q)_n}=\sum_{n=0}^\infty \frac{\Phi_n\left(\lambda, \lambda_1\right)\Phi_n\left(\lambda_2, \lambda_3, \lambda_4\right) }{(q;q)_n}.
\end{align}
\end{tl}

 In the same vein as above, we also reformulate \eqref{zzz-uuu} as a new symmetric expansion for ${}_3\phi_2$ series. It makes clear that ${}_3\phi_2$ series has an automorphism subgroup
 $S_2\times S_2\times S_2$, not only the automorphism group $S_5$.
\begin{dl}
\begin{align}     &  { }_3 \phi_2\left[\begin{array}{c}
\lambda_2 \lambda_3, \lambda_2 \lambda_4, \lambda_3 \lambda_4 \\
\lambda_1 \lambda_2 \lambda_3 \lambda_4, \lambda_2 \lambda_3 \lambda_4 \lambda_5
\end{array} ; q, \lambda_1 \lambda_5\right]= \frac{\poq{\lambda_1 \lambda_4, \lambda_2 \lambda_4, \lambda_2 \lambda_5}{\infty}}{\poq{\lambda_1 \lambda_2 \lambda_3 \lambda_4, \lambda_2 \lambda_3 \lambda_4 \lambda_5}{\infty}}\nonumber\\
&\qquad\qquad\times \sum_{n=0}^\infty\frac{\Phi_{n}(\lambda_1 \lambda_2 \lambda_3;\lambda_1 , \lambda_2  )\Phi_{n}( \lambda_3 \lambda_4\lambda_5;  \lambda_4,\lambda_5)}{\poq{q}{n}}.
         \label{zzz-uuu-kkk}
    \end{align}
\end{dl}
\pf
It only needs to make the simultaneous replacements in \eqref{zzz-uuu}
$$
\left(x,y,\alpha,\lambda_3,\beta\right)\to\bigg(\frac{x_1 x_5}{\lambda_1},\frac{x_2
   x_5}{\lambda_1}, \frac{x_1 x_2 x_3
   x_5}{\lambda_1},\frac{\lambda_1
   x_4}{x_5},\lambda_1 x_3 x_4\bigg).
$$
The result is as follows
\begin{align*}     &  { }_3 \phi_2\left[\begin{array}{c}
x_2 x_3, x_2 x_4, x_3 x_4 \\
x_1 x_2 x_3 x_4, x_2 x_3 x_4 x_5
\end{array} ; q, x_1 x_5\right]= \frac{\poq{x_1 x_4, x_2 x_4, x_2 x_5}{\infty}}{\poq{x_1 x_2 x_3 x_4, x_2 x_3 x_4 x_5}{\infty}}\\
&\times \sum_{n=0}^\infty\frac{\Phi_{n}(x_1 x_2 x_3 x_5/\lambda_1;x_1 x_5/\lambda_1, x_2 x_5/\lambda_1 )\Phi_{n}( \lambda_1x_3 x_4; \lambda_1, \lambda_1 x_4/x_5)}{\poq{q}{n}}.  \end{align*}
Using the homogeneous property \eqref{homo} to simplify this identity and relabel $x_i$ with $\lambda_i$, we obtain the desired result.
\qed

As above, as a direct application of \eqref{zzz-uuu-kkk}, we easily find the transformation among two bilinear GFs for the ${}_3\phi_2$ series.
\begin{tl}
  \begin{align}
  \left(\lambda_1 \lambda_3, \lambda_2 \lambda_3\right)_{\infty}&\sum_{n=0}^\infty \frac{\Phi_n\left(\lambda_1, \lambda_2\right) \Phi_n\left(\lambda_3 \lambda_4 \lambda_5 ; \lambda_3, \lambda_4, \lambda_5\right)}{\poq{q}{n}}\label{VVVVVV}\\
&=\sum_{n=0}^\infty \frac{\Phi_{n}(\lambda_1 \lambda_2 \lambda_3;\lambda_1 , \lambda_2  )\Phi_{n}( \lambda_3 \lambda_4\lambda_5;  \lambda_4,\lambda_5)}{\poq{q}{n}}.\nonumber
  \end{align}
\end{tl}

\pf Recall that \eqref{yyy-yyy} states
\begin{align*}
&\sum_{n=0}^\infty\frac{\Phi_{n}(\alpha ;x, y)\Phi_{n}(\beta ;\lambda_1, \lambda_2,\lambda_3)}{\poq{q}{n}}=\frac{\poq{\alpha \lambda_3,\beta y}{\infty}}{\poq{x \lambda_1,x\lambda_3, y \lambda_2, y\lambda_3}{\infty}}\\
&\qquad\times\sum_{n=0}^\infty\frac{\poq{\alpha/x,\beta/\lambda_2,y\lambda_3}{n}}{\poq{q,\alpha \lambda_3,\beta y}{n}}\,{}_2\phi_1\left[\begin{matrix}
        \alpha/y,x \lambda_3\\\alpha \lambda_3q^n   \end{matrix};q,y\lambda_1q^n\right](x\lambda_2)^n.\nonumber
\end{align*}
Set $\alpha=0,\beta=\lambda_1\lambda_2\lambda_3$. Then it becomes
\begin{align*}        &\sum_{n=0}^\infty\frac{\Phi_{n}(x, y)\Phi_{n}(\lambda_1\lambda_2\lambda_3;\lambda_1, \lambda_2,\lambda_3)}{\poq{q}{n}}=\frac{\poq{y\lambda_1\lambda_2\lambda_3 }{\infty}}{\poq{x \lambda_1,x\lambda_3, y \lambda_2, y\lambda_3}{\infty}}\\&\qquad\times\sum_{n=0}^\infty\frac{\poq{\lambda_1\lambda_3,y\lambda_3}{n}}{\poq{q,y\lambda_1\lambda_2\lambda_3 }{n}}\,{}_1\phi_0\left[\begin{matrix}
        x \lambda_3\\-  \end{matrix};q,y\lambda_1q^n\right](x\lambda_2)^n\\
        &=\frac{\poq{y\lambda_1\lambda_2\lambda_3 ,xy\lambda_1\lambda_3}{\infty}}{\poq{x \lambda_1,x\lambda_3,y\lambda_1, y \lambda_2, y\lambda_3}{\infty}}{}_3\phi_2\left[\begin{matrix}
       \lambda_1\lambda_3,y\lambda_3,y\lambda_1\\ y\lambda_1\lambda_2\lambda_3 , xy\lambda_1\lambda_3 \end{matrix};q,x\lambda_2\right].
\end{align*}
It is equivalent to
\begin{align*}    {}_3\phi_2\left[\begin{matrix}
       \lambda_1\lambda_3,y\lambda_3,y\lambda_1\\ y\lambda_1\lambda_2\lambda_3 , xy\lambda_1\lambda_3 \end{matrix};q,x\lambda_2\right]&=\frac{\poq{x \lambda_1,x\lambda_3,y\lambda_1, y \lambda_2, y\lambda_3}{\infty}}{\poq{y\lambda_1\lambda_2\lambda_3 ,xy\lambda_1\lambda_3}{\infty}}\\
       &\times\sum_{n=0}^\infty\frac{\Phi_{n}(x, y)\Phi_{n}(\lambda_1\lambda_2\lambda_3;\lambda_1, \lambda_2,\lambda_3)}{\poq{q}{n}}.
\end{align*}Under the replacements
$(x,y,\lambda_1,\lambda_2,\lambda_3)\to(\lambda_1,\lambda_2,\lambda_3,\lambda_5,\lambda_4)$, it takes the form
\begin{align*}  & {}_3\phi_2\left[\begin{array}{c}\lambda_2 \lambda_3, \lambda_2 \lambda_4, \lambda_3 \lambda_4 \\ \lambda_1 \lambda_2 \lambda_3 \lambda_4, \lambda_2 \lambda_3 \lambda_4 \lambda_5\end{array} ; q,  \lambda_1 \lambda_5\right] \\
& =\frac{\left(\lambda_1 \lambda_3, \lambda_1 \lambda_4, \lambda_2 \lambda_3, \lambda_2 \lambda_4, \lambda_2 \lambda_5\right)_{\infty}}{\left(\lambda_1 \lambda_2 \lambda_3 \lambda_4, \lambda_2 \lambda_3 \lambda_4 \lambda_5\right)_{\infty}}\sum_{n=0}^\infty \frac{\Phi_n\left(\lambda_1, \lambda_2\right)\Phi_n\left(\lambda_3 \lambda_4 \lambda_5 ; \lambda_3, \lambda_4, \lambda_5\right)}{(q;q)_n}.\nonumber\end{align*}
A quick comparison of this with \eqref{zzz-uuu-kkk}   gives rise to \eqref{VVVVVV}.\qed

We end this paper by pointing out that Theorem
\ref{A-S-Cth1.1-genal} including its corollaries  \eqref{yyy-yyy-gel} and \eqref{yyy-yyy} deserve our great attention. Due to limitation of space, we won't go into more discussion here and postpone for forthcoming  study.

\vspace{18pt}

\noindent{\bf Funding} \qquad This research was supported by the Natural Science Foundation of Zhejiang Province under
Grant No. LY24A010012 and the National Natural Science Foundation of China under Grant Nos. 12001492 and 12471315.

\noindent{\bf Availability of data and materials}\qquad No data or material were used in producing the results in this theoretical work.

\noindent{\bf Declarations Conflict of interest}\qquad The authors declare that they have no competing interests relating to this work.

\bibliographystyle{amsplain}

\end{document}